\documentclass[11pt]{article}
\usepackage{mathrsfs}
\usepackage{latexsym,lineno}
\usepackage{epsfig}
\usepackage{color}
\usepackage{amsmath}\usepackage{fleqn}\usepackage{verbatim}\usepackage{epsf}
\usepackage{amsthm}\usepackage{graphicx, float}\usepackage{graphicx}
\usepackage{amsfonts}\usepackage{amssymb}\usepackage{graphpap}
\usepackage{epic}\usepackage{curves}
\usepackage{tikz}
\usepackage{subfigure}

\newtheorem{thm}{Theorem}

\newtheorem{prob}{Problem}

\newtheorem{remark}{Remark}

\long\def\invis#1{}
\title{\bf Solution to a Forcible Version of a Graphic \\Sequence Problem\thanks{Supported in part by the National Natural Science Foundation of China (No. 11871329)}}
\author{Mao-cheng Cai
\\Academy of Mathematics and Systems Science, \\Chinese Academy of Sciences\\
Beijing 100190, P. R. China\\ and\\
Liying Kang\thanks{Corresponding authors. Email address: lykang@shu.edu.cn (L. Kang)}\\Department of Mathematics\\Shanghai University, Shanghai 200444, P.R. China}
\begin{document}
\date{}
\maketitle
\thispagestyle{empty}
\begin{abstract}
Let $A_n=(a_1,a_2,\ldots,a_n)$ and $B_n=(b_1,b_2,\ldots,b_n)$ be nonnegative integer sequences with $A_n\le B_n$. The purpose of this note is to give a good characterization such that every integer sequence $\pi=(d_1,d_2,\ldots d_n)$ with even sum and $A_n\le \pi\le B_n$ is graphic. This solves a forcible version of problem posed by Niessen and generalizes the Erd\H{o}s--Gallai theorem.\vspace{2mm}

\vspace{1mm}
\end{abstract}

First let us introduce some terminology and notations.

Let $A_n=(a_1,a_2,\ldots,a_n)$ and $B_n=(b_1,b_2,\ldots,b_n)$ be nonnegative integer sequences with $a_i\le b_i,\,1\le i\le n$, written as $A_n\le B_n$.
A nonnegative integer sequence $\pi=(d_1,d_2, \ldots, d_n)$ is called {\em graphic} if there is some simple graph having degree sequence $\pi$.

For simplicity, let $\mathcal{S}[A_n,B_n] $ denote the set of integer sequences $\pi=(d_1,d_2,\ldots, d_n)$ with even sum and $A_n\le \pi\le B_n$.

The following Erd\H{o}s--Gallai theorem gave a good
characterization for a nonnegative integer sequence to be graphic.

\begin{thm}[Erd\H{o}s--Gallai~\cite{EG}]\label{thm-EG} \ Let
$\pi=(d_1,d_2,\cdots,d_n)$ be a nonnegative integer sequence in non-increasing order. Then $\pi$ is graphic if and only if the sum of $\pi$ is even and
\begin{equation}
\sum_{i=1}^t d_i\leq t(t-1)+\sum_{i=t+1}^n\min\{t,\,d_i\} \
\mbox{ for every $\,t$, $1\le t\le n$}.\label{L1}
\end{equation}
\end{thm}

Motivated by this theorem, Niessen posed the following
\begin{prob}{\rm(\cite{Nie})}\label{L-Nie} \
Let $A_n$ and $B_n$ be integer sequences with $0\le A_n\le B_n$.
Give a simple characterization {\rm (}\/like the
above theorem\/{\rm)} for the existence of a graphic sequence $\pi=(d_1,d_2,\ldots,d_n)\in \mathcal{S}[A_n,B_n]$.
\end{prob}

The problem is regarded as the {\em potential} version. A {\em forcible} version of the problem is the following

\begin{prob}{\rm(\cite{GY})}\label{L-GY} \ Let $A_n$ and $B_n$ be integer sequences with $0\le A_n\le B_n$.
Give a simple characterization {\rm (}\/like the
above theorem\/{\rm)} such that every sequence $\pi=(d_1,d_2,\ldots,d_n)\in \mathcal{S}[A_n,B_n]$ is graphic.
\end{prob}

For convenience, we say that $A_n$ and $B_n$ are in {\em good} order $\mathcal{A}$ (respectively, $\mathcal{B}$) if $a_i> a_{i+1}$ or $a_i=a_{i+1}$  and $b_i\ge b_{i+1}$ (respectively, $a_i\ge a_{i+1}$ and $a_i+b_i\ge a_{i+1}+b_{i+1}$) for $i=1, 2,\ldots,n-1$.

Given $A_n$ and $B_n$ in good order $\mathcal{A}$, define for $t=0,1,\ldots, n$
\begin{align*}
J(t) &= \{i\mid i\ge t+1\, , \,b_i\ge t+1\},\\
\alpha(t) &=
\begin{cases}1 & \mbox{if $a_i=b_i\,\forall\,i\in J(t)$ and $\sum\limits_{i\in J(t)}b_i+t\/|J(t)|\equiv
1\pmod{2}$},\\
0 & \mbox{otherwise.}
\end{cases}
\end{align*}

Cai et al.\,\cite{CDZ} gave a solution to Problem \ref{L-Nie}, very similar in form to Theorem \ref{thm-EG}.
\begin{thm}{\rm(\cite{CDZ})}\label{thm-CDZ} \ Let $A_n$ and $B_n$ be in good order $\mathcal{A}$.
Then there exists a graphic sequence $\pi\in \mathcal{S}[A_n,B_n]$ if and only if
\begin{equation}
\sum_{i=1}^t a_i\leq
t(t-1)+\sum_{i=t+1}^n\min\{t,\,b_i\}-\alpha(t) \ \mbox{ for
every $\,t$, $0\le t\le n$}.\label{L2}
\end{equation}
\end{thm}

Possibly inspired by a result of Niessen~\cite{Nie1} , Guo and Yin~\cite{GY} posed and studied Problem~\ref{L-GY}, obtained imperfect results for the case $A_n$ and $B_n$ in good order $\mathcal{B}$.

Given $A_n$ and $B_n$ in good order $\mathcal{B}$, define for $t=0,1,\ldots, n$
\begin{align*}
J(t) &= \{i\mid i\ge t+1\, , \,b_i\ge t+1\},\\
\xi(t) &=
\begin{cases}1 & \mbox{if $a_i<b_i$ for some $i\in J(t)$ or $\sum\limits_{i\in J(t)}b_i+t\/|J(t)|\equiv
1\pmod{2}$},\\
0 & \mbox{otherwise.}
\end{cases}
\end{align*}
\begin{thm}{\rm(\cite{GY})} \ Let $A_n$ and $B_n$ be in good order $\mathcal{B}$.
If every sequence $\pi\in \mathcal{S}[A_n,B_n]$ is graphic, then for $t=0,1,\ldots,n$,
\begin{equation}
\sum_{i=1}^t b_i\leq
\begin{cases}
t(t-1)+\sum\limits_{i=t+1}^n\min\{t,\,a_i\}-\xi(t)+2 \ & \mbox{ if $a_i< b_i$ for some $i$},\\
t(t-1)+\sum\limits_{i=t+1}^n\min\{t,\,a_i\}-\xi(t) &\mbox{ if $a_i=b_i$ for each $i$}.
\end{cases}\label{L3}
\end{equation}
\end{thm}

\begin{thm}{\rm(\cite{GY})} \ Let $A_n$ and
$B_n$ be in good order $\mathcal{B}$. If for $t=0,1,\ldots,n$,
\begin{equation}
\sum_{i=1}^t b_i\leq
\begin{cases}
t(t-1)+\sum_{i=t+1}^n\min\{t,\,a_i\}-\xi(t)+1 & \mbox{if $a_i<b_i$ for some $i$},\\
t(t-1)+\sum_{i=t+1}^n\min\{t,\,a_i\}-\xi(t) & \mbox{if $a_i=b_i$ for each $i$},
\end{cases}
\label{L4}
\end{equation}
then every sequence $\pi\in \mathcal{S}[A_n,B_n]$ is graphic.
\end{thm}

Clearly, there is a gap between the necessary and sufficient conditions given above.

In \cite{Cai-Kang} we eliminated the gap and characterized the case $A_n$ and $B_n$ in good order $\mathcal{B}$ by Theorem~\ref{thm-cai-kang}.

Given $A_n$ and $B_n$ in good order $\mathcal{B}$, define for $t=1,2,\ldots,n$
\begin{align*}
J'(t)&=\{i> t\mid a_i\ge t\},\\
\beta'(t)&=\begin{cases}
1 & \mbox{if $A_n\neq B_n$, $a_i=b_i\ \forall i\in J'(t)$}\
  \mbox{and $\sum\limits_{i=1}^t b_i+\sum\limits_{i=t+1}^n a_i \equiv 1 \pmod{2}$,}\\
0 & \mbox{otherwise.}
\end{cases}
\end{align*}

\begin{thm}{\rm(\cite{Cai-Kang})}\label{thm-cai-kang} \ Let $A_n$ and $B_n$ be in good order $\mathcal{B}$. Every sequence $\pi\in \mathcal{S}[A_n,B_n]$ is graphic if and only if
\begin{equation}
\sum_{i=1}^t b_i\leq
t(t-1)+\sum_{i=t+1}^n\min\{t,\,a_i\}+\beta'(t)\quad \mbox{for every $t$, $1\le t\le n$.}
\label{L-thm}
\end{equation}
\end{thm}

Now it should be pointed out that in good order $\mathcal{A}$ and in good order $\mathcal{B}$ are essentially different.  Given nonnegative integer sequences $A_n$ and $B_n$ with $A_n\le B_n$, it is always possible to arrange them in good order $\mathcal{A}$. But it is less likely to arrange them in good order $\mathcal{B}$ because, generally speaking, the conditions $b_1\ge b_2\ge\cdots\ge b_n$ and $a_1+b_1\ge a_2+b_2\ge\cdots\ge a_n+b_n$ are not necessarily compatible.

Therefore, Problem~\ref{L-Nie} was solved completely, but Problem~\ref{L-GY} is not, solved only for the special case $A_n$ and $B_n$ in good order $\mathcal{B}$ by Theorem~\ref{thm-cai-kang}. However, the approach used in~\cite{Cai-Kang} can be modified to deal with the general case.

The purpose of this note is to give a solution to Problem~\ref{L-GY}, similar in form to Theorem~\ref{thm-EG}.

Let $t$ be an integer with $1\le t\le n$. We say that $A_n$ and $B_n$ are in {\em good}\/ order $O(t)$ if
\begin{itemize}
\item $b_i+\min\{t,a_i\}> b_{i+1}+\min\{t,a_{i+1}\}$ or
\item $b_i> b_{i+1}$ when $b_i+\min\{t,a_i\}= b_{i+1}+\min\{t,a_{i+1}\}$ or
\item $b_i+a_i\ge b_{i+1}+a_{i+1}$ when $b_i+\min\{t,a_i\}= b_{i+1}+\min\{t,a_{i+1}\}$ and $b_i= b_{i+1}$
\end{itemize}
for $i=1,2,\ldots,n-1$.

Obviously, for each $t=1,2,\ldots,n$, $A_n$ and $B_n$ can be arranged as $A_{tn}=(a_{t1},a_{t2},\ldots,a_{tn})$ and $B_{tn}=(b_{t1},b_{t2},\ldots,b_{tn})$ such that $A_{tn}$ and $B_{tn}$ are in good order $O(t)$. We define
\begin{align*}
\rho(t)&= b_{tt}+\min\{t,a_{tt}\}, \quad J^*(t)=\{i\mid b_{ti}+\min\{t,a_{ti}\}= \rho(t)\},\\
I_1(t)&=\{1,2,\ldots,t\},\quad I_2(t)=\{i> t\mid a_{ti}\ge t\},\quad
I_3(t)=\{i>t\mid a_{ti}< t\},\\
\beta(t)&=\begin{cases}
1 & \mbox{if $A_n\neq B_n$, $a_{ti}=b_{ti}\; \forall\, i\in I_2(t)$},\
  \mbox{$\sum\limits_{i=1}^t b_{ti}+\sum\limits_{i=t+1}^n a_{ti} \equiv 1 \pmod{2}$}\\
  &\mbox{and $b_{ti}+a_{ti}\equiv 0\pmod{2}\ \forall i\in I_1(t)\cap J^*(t)$ when $I_2(t)\cap J^*(t)\neq \emptyset$,}\\
0 & \mbox{otherwise.}
\end{cases}
\end{align*}

Now let us show
\begin{equation}
b_{ti} \ge \begin{cases}
\min\{a_{tj}+1,b_{tj}\}\ge a_{tj} & \mbox{ if $i<j$,}\\
a_{tj}  &\mbox{ if $i,j\in J^*(t)$}.
\end{cases}\label{N_1}
\end{equation}
Indeed, assuming $b_{ti}<\min\{a_{tj}+1,b_{tj}\}$, then $b_{tj}>b_{ti}$, $a_{tj}\ge b_{ti}\ge a_{ti}$,
$b_{tj}+\min\{t,a_{tj}\}>b_{ti}+\min\{t,a_{ti}\}$, thus $j<i$, a contradiction. Similarly, assuming $b_{ti}<a_{tj}$, then $b_{tj}+\min\{t,a_{tj}\}>b_{ti}+\min\{t,a_{ti}\}$ but $b_{tj}+\min\{t,a_{tj}\}=b_{ti}+\min\{t,a_{ti}\}$ because $i,j\in J^*(T)$.

\begin{thm}\label{thm-cai1} Let $A_n$ and $B_n$ be integer sequences with $0\le A_n\le B_n$. Every sequence $\pi\in \mathcal{S}[A_n,B_n]$ is graphic if and only if
\begin{equation}
\sum_{i=1}^t b_{ti}\leq
t(t-1)+\sum_{i=t+1}^n\min\{t,\,a_{ti}\}+\beta(t)\quad \mbox{for every $t$, $1\le t\le n$.}
\label{L-thm1}
\end{equation}
\end{thm}
\noindent{\bf Proof}.\ We may first assume that $\mathcal{S}[A_n,B_n]\ne \emptyset$, for otherwise the theorem holds trivially. We may further assume that $A_n\ne B_n$, for otherwise $\beta(t)=0,\,t=1,2\ldots,n$, so that (\ref{L-thm1}) becomes (\ref{L1}).

{\bf Necessity}. For each fixed $t$ with $1\le t\le n$, consider an integer sequence $\pi^*=(d^*_1,d^*_2,\ldots,d^*_n)$ satisfying
\begin{equation}
\begin{cases}
d_i^* = b_{ti}  &\mbox{if $i\in I_1(t)$},\\
a_{ti}\le d_i^*\le \min\{a_{ti}+1,b_{ti}\} &\mbox{if $i\in I_2(t)$,}\\
d^*_i = a_{ti} & \mbox{if $i\in I_3(t)$}.
\end{cases}\label{L-proof}
\end{equation}
Then it follows from (\ref{N_1}) that
\begin{equation}
d_i^*\ge d_j^*\quad \mbox{for $1\le i\le t<j\le n$.}\label{L-proof1}
\end{equation}

Now we distinguish two cases.

\noindent{\bf Case 1:} There is a graphic sequence $\pi^*=(d_1^*,d_2^*,\ldots,d_n^*)\in \mathcal{S}[A_n,B_n]$ satisfying (\ref{L-proof}).

If necessary, we order $d_1^*,d_2^*,\ldots,d_n^*$ such that $d^*_{i_1}\ge d^*_{i_2}\ge\cdots\ge d^*_{i_n}$ with the result that $\{i_1,i_2,\ldots,i_t\}=I_1(t)$ in view of (\ref{L-proof1}). Since $\pi^*$ is graphic, applying Theorem~\ref{thm-EG} to $(d^*_{i_1}, d^*_{i_2},\cdots, d^*_{i_n})$, we have
\begin{align*}
\sum\limits_{i=1}^t b_{ti}=\sum\limits_{j=1}^t d^*_{i_j}&\le t(t-1)+\sum\limits_{j=t+1}^n\min\{t,d^*_{i_j}\}\\
&=t(t-1)+t|I_2(t)|+\sum\limits_{i\in I_3(t)}a_{ti}
=t(t-1)+\sum\limits^n_{i=t+1}\min\{t, a_{ti}\}.
\end{align*}

Moreover, $\beta(t)=0$ in that if $a_{ti}=b_{ti}$ for all $i\in I_2(t)$, then $\sum\limits_{i=1}^t b_{ti}+\sum\limits_{i=t+1}^n a_{ti} =\sum\limits^n_{i=1}d^*_i \equiv 0 \pmod{2}$.
\vspace*{2mm}

\noindent{\bf Case 2:} There is no such sequence $\pi^*=(d^*_1,d^*_2,\ldots,d^*_n)\in \mathcal{S}[A_n,B_n]$.

Then $a_{ti} =b_{ti}\; \forall\, i\in I_2(t)$ and $
\sum\limits^n_{i=1}d^*_i= \sum\limits_{i=1}^t b_{ti}+\sum\limits_{i=t+1}^n a_{tj} \equiv 1 \pmod{2}$,
or else Case 1 would occur. There are two subcases.

\noindent{\bf Subcase 2.1:} There are $j'\in I_2(t)\cap J^*(t)$ and $i'\in I_1(t)\cap J^*(t)$ such that $b_{ti'}+a_{ti'}\equiv 1\pmod{2}$. Then $\beta(t)=0$.

Clearly $b_{ti'}>a_{ti'}$, $b_{tj'}=a_{tj'}\ge t$ as $j'\in I_2(t)$. Thus $b_{tj'}+t=b_{tj'}+\min\{t,a_{tj'}\}=b_{ti'}+\min\{t,a_{ti'}\}$ since $i',j'\in J^*(t)$.
Then $a_{ti'}< t$ otherwise $a_{ti'}\ge t$, $b_{ti'}=b_{tj'}$, $b_{ti'}+a_{ti'}<2b_{ti'}= b_{tj'}+a_{tj'}$, yielding $j'<i'$, a contradiction. Hence
\begin{equation}
b_{tj'}+t=b_{tj'}+\min\{t,a_{tj'}\}=b_{ti'}+\min\{t,a_{ti'}\}=b_{ti'}+a_{ti'}.\label{N-1}
\end{equation}

Replace $d^*_{i'}$ and $d^*_{j'}$ in $\pi^*$ with $d^*_{j'}$ and $a_{ti'}$, respectively, and denote the new sequence by $\bar{\pi}^*=(\bar{d}^*_1,\bar{d}^*_2,\ldots,\bar{d}^*_n)$. Let us show that
\begin{equation}
\bar{d}^*_i\ge \bar{d}^*_j\quad \mbox{for $1\le i\le t<j\le n$}.\label{N-2}
\end{equation}
By (\ref{L-proof1}), (\ref{N-2}) holds if $i\neq i'$ and $j\neq j'$. As $i',j'\in J^*(t)$, then $k\in J^*(t)$ for every $k$ with $i'\le k\le j'$. Thus for $i=i'$ or $j=j'$, (\ref{N-2}) drives easily from (\ref{N_1}).

Moreover, $\sum\limits^n_{i=1}\bar{d}^*_i= \sum\limits^n_{i=1}{d}^*_i-b_{ti'}+a_{ti'}\equiv\sum\limits^n_{i=1}{d}^*_i+1\equiv 0 \pmod{2}$, thus $\bar{\pi}^*$ is graphic. Applying a similar argument used in Case 1 to  $\bar{\pi}^*$, we obtain
\begin{equation}
\sum\limits_{i=1}^t \bar{d}^*_{i}\le t(t-1)+\sum\limits_{i=t+1}^n\min\{t,\bar{d}^*_{i}\}.\label{N-3}
\end{equation}
On the other hand,
\begin{align*}
\sum\limits_{i=1}^t \bar{d}^*_{i}=\sum\limits_{i=1}^t b_{ti}-b_{ti'}+b_{tj'} \mbox{ and } \sum\limits_{i=t+1}^n\min\{t,\bar{d}^*_{ti}\}=\sum\limits_{i=t+1}^n\min\{t,a_{ti}\}-t+a_{ti'},
\end{align*}
combined with (\ref{N-1}) and (\ref{N-3}), we have
$$
\sum\limits_{i=1}^t b_{ti}\le t(t-1)+\sum\limits_{i=t+1}^n\min\{t,a_{ti}\}.
$$

\noindent{\bf Subcase 2.2:} $b_{ti}+a_{ti}\equiv 0\pmod{2}\; \forall\, i\in I_1(t)\cap J^*(t)$ provided $I_2(t)\cap J^*(t)\neq \emptyset$. Then $\beta(t)=1$.

Since $\mathcal{S}[A_n,B_n]\neq\emptyset$, there exists $i^*\in I_3(t)$ or $i^*\in I_1(t)$ such that $a_{ti^*}<b_{ti^*}$.
Replace $d^*_{i^*}$ in $\pi^*$ with $d^*_{i^*}+1$ or $d^*_{i^*}-1$ according to whether or not there exists an $i^*\in I_3(t)$ with $a_{ti^*}<b_{ti^*}$, and denote the new sequence by $\hat{\pi}^*=(\hat{d}^*_1,\hat{d}^*_2,\ldots,\hat{d}^*_n)$. Clearly $\hat{\pi}^*\in \mathcal{S}[A_n,B_n]$ as the sum of $\hat{\pi}^*$ is even, hence is graphic. Let us show that
\begin{equation}
\begin{array}{ll}
\hat{d}^*_j\ge \hat{d}^*_{i^*}  \ & \mbox{for every $j\le t$ if $i^*\in I_3(t)$,}\\
\hat{d}^*_{i^*} \ge \hat{d}^*_j \ & \mbox{for every $j> t$ if $i^*\in I_1(t)$}.
\end{array}\label{L-prf}
\end{equation}
In the case $i^*\in I_3(t)$ and $j\le t$, then $\hat{d}^*_j=b_{tj}\ge \min\{a_{ti^*}+1, b_{ti^*}\}=\hat{d}^*_{i^*}$ by (\ref{N_1}). And in the other case, $i^*\in I_1(t)$ and $a_{tj}=b_{tj}$ for every $j>t$, then $b_{ti^*}>b_{tj}$, for otherwise $a_{ti^*}<b_{ti^*}\le b_{tj}=a_{tj}$, implying $j<i^*$, a contradiction. Hence $\hat{d}^*_{i^*}=d^*_{i^*}-1\ge d^*_j=\hat{d}^*_j$.

Similarly, we order $\hat{d}^*_1,\hat{d}^*_2,\ldots,\hat{d}^*_n$ such that $\hat{d}^*_{i_1}\ge \hat{d}^*_{i_2}\ge\cdots\ge \hat{d}^*_{i_n}$, with the result that $\{i_1,i_2,\ldots,i_t\}=I_1(t)$ due to (\ref{L-prf}). Since $\hat{\pi}$ is graphic, applying Theorem~\ref{thm-EG} to $(\hat{d}^*_{i_1}, \hat{d}^*_{i_2},\cdots, \hat{d}^*_{i_n})$, we have
in the case $i^*\in I_3(t)$
\begin{align*}
\sum\limits_{i=1}^t b_{ti} & =\sum\limits_{j=1}^t \hat{d}^*_{i_j}\le t(t-1)+\sum\limits_{j=t+1}^n\min\{t,\hat{d}^*_{i_j}\}\\
& =t(t-1)+\sum\limits^n_{i=t+1}\min\{t, d^*_i\}+1=t(t-1)+\sum\limits^n_{i=t+1}\min\{t, a_{ti}\}+1
\end{align*}
and in the other case
\begin{align*}
\sum\limits_{i=1}^t b_{ti}-1 & =\sum\limits_{j=1}^t \hat{d}^*_{i_j}\le t(t-1)+\sum\limits_{j=t+1}^n\min\{t,\hat{d}^*_{i_j}\}\\
& =t(t-1)+\sum\limits^n_{i=t+1}\min\{t, d^*_i\}=t(t-1)+\sum\limits^n_{i=t+1}\min\{t, a_{ti}\}.
\end{align*}
Therefore (\ref{L-thm1}) holds in both cases.\vspace{1em}

{\bf Sufficiency}. Taking any sequence $S_n=(d_1,d_2,\ldots,d_n)\in \mathcal{S}[A_n,B_n]$, we
order $S_n$ as $d_{i_1}\ge d_{i_2}\ge\ldots\ge d_{i_n}$.

According to Theorem~\ref{thm-EG}, we need to show that
\begin{equation}
t(t-1)+\sum^n_{j=t+1}\min\{t,d_{i_j}\}-\sum^t_{j=1} d_{i_j}\ge 0\label{L-prfs}
\end{equation}
for every $t$, $1\le t\le n$.

For simplicity, let $\delta(S_n)$ stand for the left-hand side of (\ref{L-prfs}) and set
$I^*_t=\{i_1,i_2,\ldots,i_t\}$.

For a $t$-set $I_t=\{j_1,j_2,\ldots,j_t\}\subseteq \{1,2,\ldots,n\}$ we define a set function
\[
f(I_t)=t(t-1)+\sum_{j\in \overline{I_t}}\min\{t,a_{tj}\}-\sum_{j\in I_t}b_{tj}.
\]
Obviously,
\begin{equation}
\delta(S_n)\geq f(I^*_t).\label{L-prfs1}
\end{equation}

Recall that $I_1(t)=\{1,2,\ldots,t\}$. Let us show that
\begin{equation}
f(I^*_t)\geq f(I_1(t)),\label{L-prfs2}
\end{equation}
or equivalently,
\[
\sum_{i\in I_1(t)\setminus I^*_t} [b_{ti}+\min\{t,a_{ti}\}]\ge \sum_{j\in I^*_t\setminus I_1(t)} [b_{tj}+\min\{t,a_{tj}\}].
\]
Indeed, if $i\in I_1(t)\setminus I^*_t$  and $j\in I^*_t\setminus I_1(t)$, then $i\le t< j$. As $A_{tn}$ and $B_{tn}$ are in good order $O(t)$,
\begin{equation}
 b_{ti}+\min\{t,a_{ti}\}\ge\rho(t)\ge b_{tj}+\min\{t,a_{tj}\}.\label{L-prfs3}
\end{equation}
Using (\ref{L-thm1}), we have
$$
\delta(S_n)\geq f(I^*_t)\geq f(I_1(t))\ge - \beta(t).
$$
Consequently, (\ref{L-prfs}) holds if $\beta(t)=0$ or one of (\ref{L-prfs1}) and (\ref{L-prfs2}) is strict.

To complete the proof, it suffices to show that (\ref{L-prfs1}) is strict if $\beta(t)=1$ and (\ref{L-prfs2}) holds with equality.

For the case $\beta(t)=1$, by definition, we have
\begin{eqnarray}
&&a_{ti}=b_{ti}\quad \mbox{for all $i\in I_2(t)$},\label{B-1}\\
&&\sum_{i=1}^t b_{ti}+\sum_{i=t+1}^n a_{ti}\equiv 1\pmod{2},\label{B-2}\\
&& b_{ti}+a_{ti}\equiv 0\pmod{2}\quad \mbox{for all $i\in I_1(t)\cap J^*(t)$ when $I_2(t)\cap J^*(t)\neq \emptyset$.}\label{B-3}
\end{eqnarray}
And for the case (\ref{L-prfs2}) being equality, we have equality in (\ref{L-prfs3}).
Clearly, the symmetric difference $I_1(t)\Delta I^*_t\subseteq J^*(t)$.
Our next aim is to show that
\begin{equation}
\sum_{j=1}^t b_{ti_j}+\sum_{j=t+1}^n a_{ti_j}\equiv 1\pmod{2},\label{B-4}
\end{equation}
equivalently by (\ref{B-2})
\begin{equation}
\sum_{i\in I_1(t)\Delta I^*_t} \{b_{ti}+a_{ti}\}\equiv 0\pmod{2}.\label{B-5}
\end{equation}

If there is an $i'\in I_1(t)\setminus I^*_t$ such that $a_{ti'}\ge t$, then
\begin{equation}
b_{ti'}=a_{ti'}=b_{tj}=a_{tj}\quad \forall j\in I^*_t\setminus I_1(t).\label{B-6}
\end{equation}
In fact, for every $j\in I^*_t\setminus I_1(t)$, we have $b_{ti'}\ge b_{tj}$ as $i'<j$ and $i',j\in J^*(t)$, implying $a_{tj}\ge t$ and $b_{ti'}=b_{tj}$ as $b_{ti'}+\min\{t,a_{t_i'}\}=b_{tj}+\min\{t,a_{t_j}\}$. Thus $j\in I_2(t)$, by (\ref{B-1}) $b_{tj}=a_{tj}\le a_{ti'}\le b_{ti'}$, (\ref{B-6}) holds. Then (\ref{B-5}) follows from (\ref{B-3}) and (\ref{B-6}).

So we may assume that $a_{ti}<t$ for every $i\in I_1(t)\setminus I^*_t$. If $I_2(t)\cap J^*(t)\neq \emptyset$, then
$b_{ti}+a_{ti}=\rho(t)\equiv 0\pmod{2}$ for every $i\in I_1(t)\setminus I^*_t$ by (\ref{B-3}). Moreover,  $b_{ti}+a_{ti}=\rho(t)\equiv 0\pmod{2}$ for every $i\in I_3(t)\cap I^*_t$ and therefore for every
$i\in I^*_t\setminus I_1(t)$, thus (\ref{B-5}) holds. And if $I_2(t)\cap J^*(t)=\emptyset$, then $b_{ti}+a_{ti}=\rho(t)$ for every $i\in I_1(t)\Delta I^*_t$, hence (\ref{B-5}) holds.

We are now ready to show that (\ref{L-prfs1}) holds strictly. Note that $\sum_{i=1}^n d_i\equiv 0\pmod{2}$,
it follows from (\ref{B-4}) that either $\sum_{i\in I^*_t} d_i<\sum_{i\in I^*_t}b_{ti}$ or there is an $i'\in \overline{I^*_t}$ such that $a_{ti'}< d_{i'}\le b_{ti'}$. And for the latter case we claim further $a_{ti'}< t$ for otherwise $i'\in I_1(t)\setminus I^*_t$ as $i'\notin I_2(t)\cup I_3(t)$, contradicting (\ref{B-6}). Therefore (\ref{L-prfs1}) is strict, as required. This completes the proof. $~~\Box$ \vskip 3mm

\begin{remark} \
Theorem~\ref{thm-cai1} gives a simple algorithm that decides whether every $\pi\in \mathcal{S}[A_n,B_n]$ is graphic in $O(n^2\log n)$ time.
\end{remark}
\begin{remark} \
As we have shown, Theorem~\ref{thm-cai1} derives from Theorem~\ref{thm-EG}. Conversely, the latter is just a special case of the former when $A_n=B_n$.
\end{remark}


\end{document}